\documentclass{amsart}
\title{Fine polar invariants of minimal singularities of surfaces}

\author{Romain Bondil}
\address{Fakult\"at f\"ur Mathematik der Ruhr-Universit\"at, Universit\"atsstr. 150, Geb. NA~2/31, 44780 Bochum, Germany. }
\email{romain.bondil@ruhr-uni-bochum.de}

\input xy
\xyoption{all}



\usepackage{amsmath,amssymb,amsthm,amsfonts,latexsym}
\usepackage{amsthm}
\usepackage[french,english]{babel}
\newtheorem{thm}{Theorem}[section]
\newtheorem{lem}[thm]{Lemma}
\newtheorem{prop}[thm]{Proposition}
\newtheorem{cor}[thm]{Corollary}

\theoremstyle{definition}
\newtheorem{defi}[thm]{Definition}
\newtheorem{Notation}[thm]{Notation}
\newtheorem{rem}[thm]{Remark}
\newtheorem{Example}[thm]{Example}
\newtheorem{Conclusion}[thm]{Conclusion}

\newcommand{\Z}{\mathbb{Z}}

\newcommand{\Q}{\mathbb{Q}}
\newcommand{\C}{\mathbb{C}}

\newcommand{\D}{\mathbb{D}}

\newcommand{\p}{\mathbb{P}}

\DeclareMathOperator{\dist}{dist}

\DeclareMathOperator{\spec}

\DeclareMathOperator{\dgre}{deg}

\DeclareMathOperator{\mult}{mult}

\DeclareMathOperator{\emdim}{emdim}

\newcommand{\Oc}{\mathcal{O}}

\newcommand{\OS}{\Oc_{S,O}}

\newcommand{\CSO}{C_{S,O}}

\newcommand{\DS}{\Delta_{S,0}}
\newcommand{\PSO}{P_{S,0}}

\newcommand{\ED}{\lfloor D \rfloor}
\newcommand {\EZ}{\lfloor Z_\Omega \rfloor}
\DeclareMathOperator{\Hilb}{Hilb}
\DeclareMathOperator{\divis}{div}

\begin{document}

\begin{abstract}
We consider the polar curves $\PSO$ arising from generic projections of a germ $(S,0)$ of complex surface singularity onto $\C^2$. Taking $(S,0)$ to be a minimal singularity of normal surface (i.e. a  rational singularity with reduced tangent cone),  we give the $\delta$-invariant of these polar curves, as well as the equisingularity-type of their  generic plane projections, which are also the discriminants of   generic projections of $(S,0)$. These two (equisingularity)-data for $\PSO$ are described in term,  on the one side of the geometry of the tangent cone of $(S,0)$ and on the other side  of  the limit-trees introduced by T. de Jong and D. van Straten for the deformation theory of these minimal singularities. These trees give a combinatorial device for the description of the polar curve which makes it  much clearer than in our previous Note on the subject.  This previous work mainly relied on a result of M.~Spivakovsky. Here we  give a geometrical proof via deformations (on the tangent cone, and what we call Scott deformations) and blow-ups, although we need Spivakovsky's  result at some point, extracting some other consequences of it along the way. 
\end{abstract}

\maketitle

\section*{Introduction}


The local polar varieties of any germ $(X,0)$ of reduced complex analytic space were introduced by L\^e D.T. and B. Teissier in~\cite{pol-loc}. In particular, the multiplicities of the {\em general}  polar varieties  are important  analytic invariants of the germ~$(X,0)$. 

 However, as also emphasised by these authors  (see also~\cite{laRabida} p.~430--431 and  \cite{Te-Segre}), there is more geometrical information to be gained on the geometry of $(X,0)$  by considering not only the multiplicity but the (e.g. Whitney)-equisingularity class of  these general polar varieties, which can be  also shown to be an analytic invariant.
 
 In this work, we will focus on the polar curves of a two-dimensional germ $(S,0)$. 
 
 Our reference on equisingularity theory for space curves will be the  m\'emoire~\cite{BGG}.
 Of course, as opposed to the case of plane curves, there is no complete set of invariants attached to a germ of   space curve describing its equisingularity class.  As a general rule, results on equisingularity beyond the case of plane curves only make sense by considering the constancy of invariants in given families.  Here we look at the family of polar curves and will consider the following invariants~:



\begin{defi}
\label{defi-equising-data}
Our equisingularity data for a germ of space curve is both~:

(i) the value of the delta invariant of the curve, 

(ii) the equisingularity class of its generic plane projection.
\end{defi}

We recall the definitions of these notions in the text (see  def.~\ref{defi-delta-invariant} and def.~\ref{def-generic-proj-curve}). The constancy of these two invariants in a family of space curves ensures Whitney conditions and actually the stronger  {\em equisaturation condition} (cf. \cite{BGG}). 

In general, this is still partial information~:  for example, another interesting invariant for space curves, namely  the semi-group of each branch, is completely independent of this equisaturation condition.



The purpose of this paper is to describe the equisingularity data in~\ref{defi-equising-data}   for the {\em the general polar curve}  of  a class of normal surface singularities called {\em minimal}.  

These  {\em minimal singularities}  were  studied in any dimension  by J.~ Koll\'ar in~\cite{Ko}. In the normal surface cases, they are also   the  {\em rational singularities with reduced fundamental cycle} and were studied by  M. ~Spivakovsky in \cite{Sp} and T. de Jong and D. van Straten in~\cite{DJ-VS}.



For these surfaces,  we  prove  that the general polar curve  is a a union of $A_n$-plane curves singularities\footnote{hence the information about the semi-group if the branches is obvious}, where the $n$'s and the contacts between these curves  can  be deduced from the resolution graph of the surface $(*)$.  This information gives in particular a complete description of part (ii)
of the data in~\ref{defi-equising-data}, i.e of the general plane projection of the polar curve, which is also the discriminant of the general projection (the coincidence of these two concepts is  a theorem, cf. section.~\ref{sec-pol-inv}).

The information on the discriminant was already given in the Note~\cite{discri} as a consequence of a result of Spivakovsky, but the statement there was clumsy.

Here we give a much nicer device that allows to read directly the  information about  of the discriminant (or the polar curve as well) from  both the information contained in the tangent cone of these singularities and the one given by a graph deduced from the resolution graph, which is precisely the  {\em limit tree} introduced by T. de Jong and D. van Straten in their study of the deformation theory for these minimal singularities (see~\cite{DJ-VS}).

We also provide an inductive proof relying much more on the geometrical properties of these minimal singularities. This proof makes the core of the paper. It still uses Spivakovsky's theorem however, mainly through a characterisation of generic polar curve on the resolution which we deduce along the way.

The several plane branches of the polar curve  lie in distinct planes in a bigger linear space, and the value of the delta invariant (part (i) in~\ref{defi-equising-data}) gives some (partial) information on the configuration of  these  planes in the space. We explain how this delta invariant is easily computed 
from the what we call the Scott deformation of the surface, which turns out to give a delta-constant deformation   of the polar curve onto bunches of generic configurations of lines.



\medskip

\noindent{\bf Organisation of the paper~:} 

\smallskip

In section~\ref{sec-pol-inv}, we recall the definition of the general  polar curve $\PSO$ of  a  germ of surface $(S,0)$, of the discriminant $\DS$ of a generic projection of $(S,0)$ onto $\C^2$ and the important result that   $\DS$ is a {\em generic} projection of the curve $\PSO$.

Section~\ref{sec-defi-mini} gives the definition of minimal singularities in general, the particular case of normal surfaces, and their characterization by their dual resolution graph. We then define in section~\ref{sec-more-dual-graph} a notion of  {\em height}  on the vertices of this resolution graph, which was used in other places such as \cite{Sp} and \cite{DJ-VS}, and corresponds to the number of  point blow-ups necessary to let  the corresponding exceptional component appear. We also give there our convention in representing dual graphs with $\bullet$ and $*$  and define {\em reduced dual graphs}  to be  the ones  in which the self-intersection for components of the tangent cone has the minimum absolute value.

In section~\ref{sec-mark}, we give the description of generic polar curves on a resolution of a minimal singularity as proved by M. Spivakovsky (thm.~\ref{thm-mark}). This result will play the following somehow different roles in the sections following it~:

\noindent (i) Section~\ref{sec-la-Note} explains how, using the full strength of this result, one may derive quite quickly a description of the generic discriminant $\DS$ (more precisely of the statement (*) for the polar curve). This sums-up  the Note~\cite{discri} in an improved way, and an mistake in an example there is corrected.

\noindent (ii) In section~\ref{sec-delta-inv}, we mention how,  using a result of Giraud,  theorem~\ref{thm-mark} also permits, at least theoretically, to deduce the $\delta$-invariant for the general polar curve from the shape of its transform on the minimal resolution.  This result is however not useful for concrete computations, for which we  use another approach in section~\ref{sec-scott}.

\noindent (iii) In section~\ref{sec-char-resol}, we get,  as a purely qualitative consequence of (i) and (ii), a characterisation of generic polar curves on  a resolution $X_N$  of both the singularity $(S,0)$ and its Nash-blow-up.  This will be the application of Spivakovsky's result we will use in the proof of our main result.

Section~\ref{sec-contrib-tg-cone} to \ref{sec-scott} form the core of the text~:

\noindent $\bullet$ in section~\ref{sec-contrib-tg-cone}, the polar curve for the tangent cone of a minimal singularity is made geometrically explicit and through the process of deformation onto the tangent cone is also seen as part of the polar curve of the singularity.

\noindent $\bullet$ in section~\ref{sec-limit-tree}, we recall what we need from  the limit tree construction of de Jong and van Straten. With this,

\noindent $\bullet$ in section~\ref{sec-discri-limit-tree}, we give,  and prove,  our main theorem describing the information in $(*)$ using the limit tree construction and the contribution of the tangent cone. 

\noindent $\bullet$ in section~\ref{sec-scott}, we show how a  special deformations of minimal singularities has a  nice interpretation in  our descriptions of polar curves and also give a nice way to compute the delta invariant of these, completing the information in  def.~\ref{defi-equising-data} (i).

This leads us to question  wether (part of) the deformation theory of these minimal singularities of surfaces could be recover  the one  of their discriminant?

\smallskip {\noindent \bf Acknowlegdement~--} The author thanks L\^e D.T. for arising the question treated here, M. Merle and  M. Spivakovsky for their remarks on~\cite{discri}, T. de Jong for pointing to us his limit-tree construction,  H. Flenner and B. Teissier for helpful conversations.

\section{Polar invariants of a normal surface singularity}
\label{sec-pol-inv}
\subsection{The general polar curve as an analytic invariant}

We  recall here the definition of the local polar variety  of a germ of surface  following~\cite{pol-loc}~:

\label{s-discri} 
Let $(S,0)$ be a normal complex surface singularity $(S,0)$, embedded in $(\C^N,0)$: for any $(N-2)$-dimensional vector subspace $D$  of $\C^N$, we consider a linear projection $\C^N \rightarrow \C^2$ with kernel $D$ and denote by 
$p_D : (S,0)\rightarrow (\C^2,0),$ the restriction of this projection to $(S,0)$.

Restricting ourselves to the $D$ such that $p_D$ is finite, and considering a small representative $S$ of the germ $(S,0)$, we define, as in \cite{pol-loc}~(2.2.2), the {\em polar curve} $C(D)$ of the germ $(S,0)$ for the direction $D$, as the closure in  $S$ of the critical locus of the restriction of $p_D$ to $S\setminus\{0\}$. It is a reduced analytic curve.

As explained in loc. cit., it makes sense to say that for an open dense subset of the Grassmann manifold $G(N-2,N)$  of  $(N-2)$-planes in $\C^N$, the space curves $C(D)$ are {\em equisingular} e.g. in   terms of Whitney-equisingularity (or strong simultaneous resolution, but this is the same for families of space curves, cf.~\cite{BGG}).  We call this equisingularity class the {\em general polar curve} for $(S,0)$ embedded in $\C^N$.

One may then compare the general polar curves obtained by two distinct embeddings of the surface into a $(\C^N,0)$ and it turns out that they are still Whitney-equisingular: this is essentially proved in \cite{laRabida} (see p.~430)  in a much more general setting (arbitrary dimension and ``relative'' polar varieties). Summing-up this~:

\begin{thm}
\label{thm-def-pol-gen}
The Whitney equisingularity-type of the general polar curve  $C(D)$ depends only on the analytic type of the germ $(S,0)$.
 \end{thm}
 
 In this paper, following somehow the program in \cite{Te-Segre}, we want to study this invariant $C(D)$ for a 
 special class of surface singularities. 
 
 \subsection{The generic discriminant as a derived  invariant}
 
  With the same notation as before, we define the {\em discriminant} $\Delta_{p_D}$ as (the germ at $0$ of) the reduced analytic curve of $(\C^2,0)$ image of the polar curve $C(D)$ by the finite morphism $p_D$. 

Again, one may show that,  for a generic choice of $D$, the discriminants obtained are {\em equisingular germs of plane curves}, and that this in turn defines an analytic invariant of $(S,0)$.

{\em We will denote  $\Delta_{S,0}$ the equisingularity class of the discriminant of a generic projection of $(S,0)$.}
 

A first advantage of $\DS$, as a germ of plane curve, is that its equisingularity class is well-defined in terms of classical invariants such as the Puiseux pairs of the branches and the intersection numbers 
between branches (cf. e.g. the introduction of \cite{BGG} for references on this subject).


As it turns outs, there is a very nice relationship  between the general polar curve and $\DS$.
For this we recall the following~:


\begin{defi}
\label{def-generic-proj-curve}
{\rm Let $(X,0)\subset (\C^N,0)$  be a germ of reduced curve. Then a linear projection $p : \C^N \rightarrow \C^2$ will said to be {\em generic}  with respect to  $(X,0)$ if the kernel of $p$ does not contain any limit of bisecants to $X$ (cf. \cite{BGG} for an explicit description of  the cone $C_5(X,0)$ formed by the limits of bisecants to $(X,0)$).}
\end{defi}

Then the equisingularity type of the germ of plane curve $(p(X),0)$ image  of $(X,0)$ by such a generic projection is uniquely defined by the saturation of the ring $\Oc_{X,0}$ (cf. \cite{BGG}).

We now  state the following transversality result (proved for curves on surfaces of $\C^3$ in \cite{BH} theorem~3.12  and in general as the ``lemme-cl\'e'' in \cite{laRabida}~V~(1.2.2)) relating polar curve and discriminants~:

\begin{thm}
\label{thm-transv-polaire}
Let $p_D~: (S,0)\rightarrow (\C^2,0)$ be as above, and $C(D)\subset (S,0)\subset (\C^N,0)$ be the corresponding polar curve. Then there is an open dense  subset $U$ of $G(N-2,N)$ such that for $D\in U$ the restriction of $p$ to $C(D)$ is generic in the sense of definition~\ref{def-generic-proj-curve}.
\end{thm}

\begin{defi}
\label{defi-PSO}
Let us  define $\PSO$ to be not the {\em Whitney-equisingularity} class of the general polar curve as in thm.~\ref{thm-def-pol-gen}, but the {\em equisaturation class} of the general polar curve (which may  be a smaller class). As we recalled after definition~\ref{defi-equising-data}, this class is precisely given by the constancy of the invariants there.
Then, the foregoing theorem~\ref{thm-transv-polaire} states that $\DS$ is {\em the generic plane projection} of  $\PSO$.
\end{defi}

As said in the introduction, the goal of this work is to determine $\PSO$ completely.






\section{Definition of minimal singularities}
\label{sec-defi-mini}
We begin with a definition valid in any dimension (following \cite{Ko} \S~3.4)~:

\begin{defi}
\label{defi-kollar}
We call a singularity  $(X,0)$ {\em minimal} if it is reduced, Cohen-Macaulay, the multiplicity  and embedding dimension of $(X,0)$ fulfil ~:

i) $\mult_0 X=\emdim_0 X -\dim_0 X +1$,

ii)  and the tangent cone $C_{X,0}$ of $X$ at $0$ is reduced.
\end{defi}

 Considering normal surfaces, one has the following characterisation~:

\begin{thm}
\label{thm-mini-equiv-rat-reduced}
Minimal singularities of normal surface  are exactly the rational surface singularities with reduced fundamental cycle (with the terminology of~\cite{Ar}).
\end{thm}

Condition (i) follows for any rational surface singularity from Artin's formulas for multiplicities and embedding dimension. Condition (ii) follows from the fact that the fundamental cycle of rational singularities is also the cycle defined by the maximal ideal. Conversely, the fact that minimal normal singularities are rational is proved in \cite{Ko} 3.4.9.  The proof that  ``reduced tangent cone'' implies 
``reduced fundamental cycle" is easy (after  our thm.~\ref{thm-Tyurina} or see e.g. in~\cite{Wahl} p.~245).

Taking $(S,0)$ to be a normal surface  singularity and  $\pi \,:\, (X,E)\rightarrow (S,0)$ to be the {\em minimal} resolution of the singularity, one associates as usual to the exceptional curve configuration 
$E=\pi^{-1}(0)$ a dual graph $\Gamma$ where each irreducible component $L_i$ in $E$ is represented by a vertex and two vertices are connected by an edge if, and only if, the corresponding 
components intersect.

Each vertex $x$ of ~$\Gamma$ (we will frequently abuse notation and write $x\in\Gamma$) is given a weight $w(x)$ defined as~:
$$w(x):=-L_x^2,$$
where $L_x^2$ is the self-intersection of the corresponding component $L_x$ on $X$.\par

For any rational surface singularity, it is well-known that all the $L_i$ are smooth rational curves and that $\Gamma$ is a tree. But in general, it takes some computation to check whether a given tree is the dual tree of a rational singularity (cf. \cite{Ar}).

On the contrary, one reads  at first sight from the dual graph that a surface singularity is minimal (cf.
\cite{Sp} II 2.3)~:
\begin{rem}
\label{rem-carac-mini-graph}
Let $\Gamma$ be any weighted graph. It is the dual graph of resolution of a minimal singularity if, and only if, $\Gamma$ is a tree and for each vertex $x\in \Gamma$ one has the following inequality~:
$$w(x)\geq v(x),$$
where $v(x)$ denotes the valence of $x$ i.e. the number of edges attached to $x$.
\end{rem} 

\section{More about the dual graphs}
\label{sec-more-dual-graph}

In the representation of the dual graph $\Gamma$ of a minimal singularities we will distinguish between the vertices with $w(x)=v(x)$ and the others.  

\begin{Notation}
\label{notation-graphe}
In representing the dual graphs of minimal singularities, we chose to represent with a $\bullet$ the vertices with $w(x)=v(x)$ so that there is no need to mention the weight above them.

On  the contrary we enumerate as $x_1,\dots,x_k$ the vertices with $w(x_i)>v(x_i)$, and let them figure as  $*$ on the graph. One should then mention the weights of the $(x_i)$ to define the graph.

In this work , we will pay much attention to the minimal singularities with the property that  for all vertices  $x_i$ with $w(x_i)>v(x_i)$ one has in fact the equality $w(x_i)=v(x_i)+1$.

Let us here call  {\em reduced} the  graphs  with this property : it is then clear that in representing these dual graphs, there is no longer need to mention the weights.
\end{Notation}

For example, to say that the  graph in figure~\ref{fig-bambou-reduit} is {\em reduced}  amounts to say that  
$w(x_1)=w(x_n)=2$ and $w(x_i)=3$ for $1<i<n$ (and the vertices with $\bullet$ all have weight two here).

\begin{figure}[h]
\caption{\label{fig-bambou-reduit}}
\setlength{\unitlength}{.7mm}
\begin{center}
\begin{picture}(0,0)(150,30) 
\put(120,22){\makebox(0,0){$* $}}
\put(120,26){\makebox(0,0){$\scriptstyle (x_1)$}}
\put(121,22){\line(1,0){8.5}}
\put(130.5,22){\circle*{2}}
\put(136,22){\makebox(0,0){$\cdots$}}
\put(140,22){\circle*{2}}
\put(141,22){\line(1,0){8,5}}
\put(150,22){\makebox(0,0){$ * $}}
\put(150,26){\makebox(0,0){$\scriptstyle (x_2)$}}

\put(151,22){\line(1,0){8.5}}
\put(160.5,22){\circle*{2}}
\put(166,22){\makebox(0,0){$\cdots$}}
\put(170,22){\circle*{2}}
\put(171,22){\line(1,0){8.5}}
\put(180,22){\makebox(0,0){$*$}}
\put(180,26){\makebox(0,0){$\scriptstyle (x_3)$}}
\put(181,22){\line(1,0){8.5}}
\put(190.5,22){\circle*{2}}
\put(196,22){\makebox(0,0){$\cdots$}}
\put(202.5,22){\makebox(0,0){$\cdots$}}
\put(207,22){\circle*{2}}
\put(208,22){\line(1,0){8.5}}
\put(217.5,22){\makebox(0,0){$ *$}}
\put(217.5,26){\makebox(0,0){$\scriptstyle (x_n)$}}
\end{picture}
\end{center}
\bigskip
\end{figure}

The geometrical meaning  of this distinction between vertices comes from the~:

\begin{thm}[Tyurina, cf.~\cite{Ty}]
\label{thm-Tyurina}
Let $(S,0)$ be a rational surface singularity and $\pi \,:\, (X,E) \rightarrow (S,0)$ its minimal resolution.
Let $b\,:\, S_1 \rightarrow S$ the blow-up of $0$ in $S$. 

Then there is a morphism $r \,:\, X \rightarrow S_1$ such that $\pi=b\circ r$ and a component $L_i$ of $E=\pi^{-1}(0)$ is contracted to a point by $r$ if, and only if,   the intersection $(L_i\cdot Z)=0$, where $Z$ is the fundamental cycle.
\end{thm}
Of course,  the components of $E$ which are not contracted by $r$ are the (strict transform by $r$) of the components of the $\p(C_{S,0})$ appearing on $S_1$.

When $(S,0)$ is a minimal singularity, the fundamental cycle is $Z=\sum_{x\in \Gamma} L_x$
and hence for a given vertex $y\in \Gamma$ the intersection $(L_y\cdot Z)$ is just $v(y)-w(y)$.

This should justify the~:
\begin{defi}
\label{defi-height}
Let $(S,0)$ be a minimal normal surface singularity and $\Gamma$ be the dual graph of its minimal resolution. 

We will say that a vertex $x$ in $\Gamma$ has {\em height one} if $w(x)>v(x)$, which from the foregoing remarks means that the corresponding component $L_x$ corresponds to a component of (the proj of) the tangent cone $C_{S,0}$.
Hence we will denote by $\Gamma _{TC}=\{x_1,\dots,,x_n\}$ the set of these vertices.

Then, we define the height of any vertex $x$ in $\Gamma$ as the number $s_x$ defined by~:
$$s_x:=\dist (x,\Gamma_{TC})+1,$$
where dist is the distance on the graph (number of edges on the geodesic between two vertices).
\end{defi}

The reader should check that this height corresponds to the number of blow-ups necessary to make the corresponding component ``appear".\footnote{This latter notion is studied more systematically  for any rational singularity as ``desingularization depth" in \cite{Le-To}; of course in this general case, it is not given  directly from a distance!}  The notation $s_x$ here comes from ~\cite{Sp} II 5.1 and was the one used in  the previous work~\cite{discri}.

\begin{figure}[h]
\caption{\label{Graphe-a-deux-limites} Minimal graph with the heights for  the vertices.}

\begin{center}
\unitlength=1cm
\begin{picture}(6,8)
\put(1,5){$*$}
\put(1,5.3){$x_1$}
\put(1,4.7){1}
\put(1.1,5.1){\line(1,0){1}}
\put(2.1,5.1){\circle*{.1}} 
\put(2.1,4.7){2}
\put(2.1,5.1){\line(1,0){1}}
\put(3.1,5.1){\circle*{.1}}
\put(2.8,4.7){3}
\put(3.1,5.1){\line(0,-1){1}}
\put(3.1,4.1){\circle*{.1}}
\put(2.8,4.1){2}
\put(3.1,4.1){\line(0,-1){1}}
\put(3,3){$*$}
\put(2.8,3){1}
\put(3.3,3){$x_2$}
\put(3.1,3.1){\line(0,-1){1}}
\put(3.1,2.1){\circle*{.1}} 
\put(2.8,2.1){2}
\put(3.1,2.1){\line(0,-1){1}}
\put(3,1){$*$}
\put(2.8,1){1}
\put(3.3,1){$x_3$}
\put(3.1,5.1){\line(1,0){1}}
\put(4.1,5.1){\circle*{.1}}
\put(4.1,4.7){2}
\put(4.1,5.1){\line(1,0){1}}
\put(5,5){$*$}
\put(5,5.3){$x_4$}
\put(5,4.7){1}
\end{picture}
\end{center}
\end{figure}
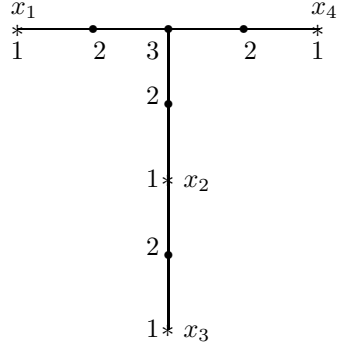

\begin{Example}
\label{Example-montre-height}
As an example, we let figure the heights on the graph in Figure~\ref{Graphe-a-deux-limites}, where the $(x_i)$ are as before the vertices  of height one (with $*$):
\end{Example}

We will also need the following~:
\begin{defi}
\label{defi-Tyurina-comp}
Let $\Gamma$ be a minimal graph. The connected components $\Gamma_i$ (for $i=1,\dots,r$) of $\Gamma\setminus \Gamma_{TC}$ are called the {\em Tyurina components} of $\Gamma$.

Theorem~\ref{thm-Tyurina} hence state that the  blow-up $S_1$ of $(S,0)$ has exactly $r$ singularities  $(S_1,O_i)$ which are minimal singularities with dual resolution graph~$\Gamma_i$.
\end{defi}

\section{A result of Spivakovsky}
\label{sec-mark}


To state this result,  we introduce a further terminology:

Let $\pi~: (X,E) \rightarrow (S,0)$ be the minimal resolution of the singularity $(S,0)$, where $E=\pi^{-1}(0)$ is the exceptional divisor, with components $L_i$. A {\em cycle} will be by definition a divisor with support on $E$ i.e. a linear combination $\sum a_i L_i$ with $a_i\in \Z$ (or $a_i\in \Q$ for a $\Q$-cycle).

Let $\Gamma$ be the dual graph of the minimal resolution $\pi$ and  for each vertex $x$ let $s_x$ denote the {\em height} defined in def.~\ref{defi-height}.

\begin{defi}
\label{defi-central}
Let then $x,y$ be two adjacent vertices on $\Gamma$~: the edge ($x,y)$ in $\Gamma$ is called a {\em central arc} if $s_x=s_y$. A vertex $x$ is called a {\em central vertex} if there are at least two vertices $y$ adjacent to $x$ such that $s_y=s_x-1$ (cf.~\cite{Sp}).

\end{defi}

We then  define a $\Q$-cycle $Z_\Omega$ on the minimal resolution $X$ of $(S,0)$ by~:

\begin{equation}
\label{eq-ZOmega}
Z_\Omega=\sum_{x\in \Gamma} s_x L_x - Z_K,
\end{equation}
\noindent where $\Gamma$ is the dual graph, and $Z_K$ is the numerically canonical $\Q$-cycle \footnote{uniquely defined by the condition that for all $x\in \Gamma,  Z_K.L_x=-2-L_x^2$ since the intersection product on $E$ is negative-definite.}.

The theorem from \cite{Sp} (theorem 5.4) is now~:
\begin{thm}
\label{thm-mark}
Let $(S,0)$ be a {\em minimal} normal surface singularity. There is a open dense subset $U'$ of the open set $U$ of theorem~\ref{thm-transv-polaire}, such that for all $D\in U'$ the strict transform $C'(D)$ of $C(D)$ on $X$~:

\noindent a) is a multi-germ of smooth curves intersecting each component $L_x$ of $E$ transversally in exactly $-Z_\Omega.L_x$ points, 

\noindent b) goes through the  point of intersection of $L_x$ and $L_y$ if and only if $s_x=s_y$ (point corresponding to a {\em central arc} of the graph). Furthermore, the curves $C'(D)$, with $D\in U'$  do not share other common points (base points) and these base points   are simple,  i.e. the curves $C'(D)$ are separated when one blows up these points once.
\end{thm}


Referring to loc. cit. for unexplained terminology, let us make the following~:
\begin{rem}
\label{rem-Nash}
Blowing-up once the base points referred to in the b) above, one gets a resolution $X_N$ of the Nash blow-up of the germ $(S,0)$. The map from $X_N$ to  the normalized Nash blow-up $N(S)$ is simply the contraction of the exceptional components  which are not intersected  by a  branch of the generic polar curve.
\end{rem}

\section{First description of the polar curve and  the discriminant}
\label{sec-la-Note}

This section essentially describes the results obtained in \cite{discri} in an improved form.
We refer to this Note for the proofs of the following lemmas~:

\begin{lem}
\label{lem-res-cp}
Let $(S,0)$ be a {\em minimal} normal surface singularity and $\pi~: X\rightarrow (S,0)$ its minimal resolution. It is known that $\pi$ is (the restriction to $S$ of) a composition  $\pi_1\circ\cdots\circ \pi_r$ of point blow-ups. Then, this composition of blow-ups  is also the  minimal resolution of the generic polar curve $C(D)$ for $D\in U'$ as in theorem~\ref{thm-mark}.
\end{lem}




The following is a slightly more precise version of loc. cit. lem.~3.2~:
\begin{lem}
\label{lem-mult-two}
For $D\in U'$ as in theorem~\ref{thm-mark}, the polar curve $(C(D),0)$ on ($S,0)$ is  a union of germs of curves of multiplicities two. In particular, it has only  smooth branches and branches of multiplicity two, the latter being exactly those for which the  strict transform goes through a central arc as in b) of theorem~\ref{thm-mark}.
\end{lem}





Let us now make a perhaps not so standard definition~:
\begin{defi}
\label{defi-contact}
Let $(\Gamma_1,0)$ und $(\Gamma_2,0)$  be two analytically irreducible curve germs  in $(\C^N,0)$. We will {\em hereafter} call {\em contact} between the $\Gamma_i$ simply the number of points blow-ups necessary to separate these two branches.
\end{defi}

For the   description of the polar curve, just recall that one calls {\em an $A_n$-curve} a curve analytically isomorphic to the plane curve defined by $x^2+y ^{n+1}=0$~:

\begin{prop}
\label{prop-descri-pol}
Let $(S,0)$ be a minimal surface singularity and $C=C(D)$ be a generic polar curve corresponding to $D$ in the open set $U'$ of thm.~\ref{thm-mark}.

Then  if $C=\cup_i \Gamma_i$ is the decomposition of $C$ into  analytic branches, denote $L_{\Gamma_i}$ the irreducible exceptional component on the minimal resolution of $S$ which intersects the strict transform of $\Gamma_i$. It is unique except in the case of central arcs. In this case just choose one between the two intersecting components. Then~:

\noindent {\rm (i)} The contact between $\Gamma_i$ and $\Gamma_j$ in the sense of def.~\ref{defi-contact} above is the minimum height in the chain between $L_{\Gamma_i}$ and $L_{\Gamma_j}$ (cf. def.~\ref{defi-height}).

\noindent {\rm (ii)} We may write rather $C$ as a union of $C=\bigcup C_i$ of curves of multiplicity two 
by taking by pairs branches intersecting the same exceptional component on $X$ that we will now denote $L_{C_i}$.

Then, each $C_i$ is a  $A_{n_i}$-curve, when the number $n_i$ equals 
$2.s(L_{C_i})$   if $C_i$ goes through a central arc, and $2.s(L_{C_i})-1$ otherwise (which comprises the case of central vertices and components of height $s$ equal to one).

We may obviously define the contact between these $A_{n_i}$-curves just by taking one branch in each, so that it is still given by (i).
\end{prop}

\begin{proof}
The statement about the contact in (i) follows from lemma~\ref{lem-res-cp}. 
The first statement in (ii) is lemma~\ref{lem-mult-two}.

Any curve of multiplicity two is a $A_n$-curve, see e.g.~\cite{BPV} p.~62. The statement about the $n_i$ follows from (i) just like the statement about the contacts.
\end{proof}


The result in prop.~\ref{prop-descri-pol} gives a complete description of the equisingularity class of the discriminant plane curve in $(\C^2,0)$ using theorem~\ref{thm-transv-polaire}~\footnote{This is  an equivalent, but more simply expressed, version of  the statement in~\cite{discri}  Cor.~4.3.}~:

\begin{prop}
\label{prop-descri-discri}
The discriminant $\Delta_{p_D}=p_D(C(D))$ has exactly the same properties as the polar curve $C(D)$ in prop.~\ref{prop-descri-pol}. This describes the generic discriminant $\DS$ as a union of $A_{n_i}$-curves with the $n_i$ and the contacts described in~\ref{prop-descri-pol}.
\end{prop}

\begin{proof}
The curves $C_i$ in  prop.~\ref{prop-descri-pol} being plane curves, they are their own generic plane projections. Hence by thm.~\ref{thm-transv-polaire}, the image $\Delta_{p_D}$ of $C(D)$ by the {\em generic projection} $p_D$ decomposes as the same union of $A_{n_i}$-curves.

We give here a direct argument to prove that the contact (in sense of  def.~\ref{defi-contact}) between the branches in $\Delta_{p_D}$ is the same as the one in $C(D)$ (in~\cite{discri}, we invoked a bilipschitz invariance which is perhaps not obvious with our definition of contact)~:
 


Considering a pair $\Gamma_1,\Gamma_2$ of branches in $C(D)$,  we may   embed $\Gamma_1\cup \Gamma_2$ into a $(\C^3,0)$ and choose coordinates so that 
$\Gamma_1$ is parametrized by $(x=t^{\varepsilon_1},y=t^n,z=0)$ and (unless the contact is one) $\Gamma_2$ is parametrized by $(x=t^{\varepsilon_2},y=0,z=t^m)$, with $\varepsilon_i=1$ or $2$.
  We then leave it to the reader that the projection defined by  $(x,y+z)$ is transverse to the $C_5$ of def.~\ref{def-generic-proj-curve} and that the contact in our sense is preserved.
\end{proof}

The foregoing  description of the discriminant still  involves the computation of the number of branches on each central vertex by Spivakovsky's formula.  We will describe a much better and condensed one in section~\ref{sec-discri-limit-tree}, which does not involve any computation and is geometrically more significant. Before, the author  would like  to make amend to the readers of 
~\cite{discri} for a mistake in the following~:

\begin{Example}[Correct version of Example 1 in~\cite{discri}]
\label{Ex-la-Note}
Consider $(S,0)$ with dual graph $\Gamma$ as on figure~\ref{fig-expl-corrige}, where, following the convention of section~\ref{sec-more-dual-graph}, the $\bullet$ denote vertices with $w(x)=v(x)$, and the others form $\Gamma_{TC}=\{x_1,\dots,x_4\}$ with the weights indicated on the graph.

The branches of the polar curve going through the components of $\Gamma_{TC}$ are just four branches going through $L_{x_1}$, which gives in the equisingularity class $\DS$ four  distinct lines through the origin with contact one with any other branch of $\DS$.

Then we have two central vertices (of height $3$ and $2$) and a central arc (with boundaries of heigth $2$), which give respectively a $A_5$ a $A_3$ and a $A_4$-curve from prop.~\ref{prop-descri-pol}  and \ref{prop-descri-discri} above.

The contact between the $A_5$ and the $A_4$ is {\em two} (and not $3$ as claimed in loc. cit.) and their contact with the other $A_3$ is one. 

Hence using coordinates, one may take as representative of the equisingularity class of $\Delta_{S,0}$ can be choosen to be~:
$$\underbrace{(x^4+y^4)}_{\mbox{The two}\; A_1}(x^2+y^5)\underbrace{(x+y^2+i.y^3).(x+y^2-iy^3)}_{\mbox{The}\; A_5}(x^2+y^5)(y^2+x^4)=0.$$
\end{Example}

\vspace{-1cm}

\begin{figure}[h]
\caption{\label{fig-expl-corrige}}

\setlength{\unitlength}{.7mm}
\begin{center}
\begin{picture}(0,0)(150,30) 

\put(48,30){\makebox(0,0){$\scriptstyle (x_1)$}}
\put(48,34){\makebox(0,0){$\scriptstyle 4$}}
\put(51.5,30){\line(1,0){6}}
 \put(57.5,30){\circle*{2}}
\put(58.5,30){\line(1,0){6}}
\put(65.5,30){\circle*{2}}
\put(66.5,30){\line(1,0){6}}
\put(76,30){\makebox(0,0){$\scriptstyle (x_2)$}}
\put(76,34){\makebox(0,0){$\scriptstyle 3$}}
\put(79.5,30){\line(1,0){6}}
\put(86.5,30){\circle*{2}}
\put(87.5,30){\line(1,0){6}}
\put(97,30){\makebox(0,0){$\scriptstyle (x_3)$}}
\put(97,34){\makebox(0,0){$\scriptstyle 2$}}
\put(65.5,29){\line(0,-1){6}}
\put(65.5,22){\circle*{2}}
\put(65.5,21){\line(0,-1){6}}
\put(65.5,16){\circle*{2}}
\put(65.5,15){\line(0,-1){6}}
\put(65.5,07.5){\makebox(0,0){$\scriptstyle (x_4)$}}
\put(71.5,07.5){\makebox(0,0){$\scriptstyle 2$}}

\end{picture}

\end{center}
\end{figure}

\vspace{1cm}

\section{The delta invariant of the polar curve}
\label{sec-delta-inv}

\begin{defi}
\label{defi-delta-invariant}
Let  $(C,0)$ be  a  germ of reduced complex curve singularity.
Let $n : \overline{C}\rightarrow C $ be its  normalisation map, which provides a finite inclusion of the local ring $\Oc_{C,0}$ into the semi-local ring $\Oc_{\overline{C}}$.

The $\delta$-invariant of $(C,0)$ is by definition $\delta (C,0):=\dim_\C  \Oc_{\overline{C}}/\Oc_{C,0}.$
\end{defi}

In the paper~\cite{Giraud}, J. Giraud gives a  way to compute $\delta(C,0)$  for any curve  on a  {\em rational} surface singularity $(S,0)$ if one knows a resolution of the surface singularity where $C'$ is a multi-germ  of smooth curves.

 To quote this result, we need the following  lemma,  proved in loc. cit. 3.6.2~:
 \begin{lem}
 \label{lem-Giraud}
 Let $p: (X,E)\rightarrow (S,0)$ be a resolution of a normal surface singularity~$(S,0)$, with $E=\pi^{-1}(0)=\cup_i E_i$. Let $D=\sum_i a_i E_i$ be  a $\Q$-cycle on $X$.
 
 There exists a unique $\Z$-cycle $V=\sum_i \alpha_i E_i$ with the property that the intersection $(V\cdot E_i)\leq  (D\cdot E_i)$ for all $i$'s, and minimum for this property.
 
 This $\Z$-cycle will be denoted as $\lfloor D\rfloor $. 
  \end{lem}
 
(In the previous lemma, minimum means that any other $\Z$-cycle with this property has the form  $\ED+ W$ with $W$ a cycle with non-negative coefficients.)
 
 In the situation of lemma~\ref{lem-Giraud}, let's associate to any curve $(C,0)\subset (S,0)$ a $\Q$-cycle $Z_C$ uniquely defined by the condition that for all irreducible component $E_i$ of $E$ the intersection number $(E_i\cdot Z_C)$ equals $(E_i\cdot C')$ where $C'$ denote the strict transform  of 
 $C$ on $X$. We may then quote (cf. loc. cit. cor.~(3.7.2))~:
 \begin{thm}
 \label{thm-Giraud}
 Let $p: (X,E) \rightarrow (S,0)$ be a resolution of a {\em rational} surface singularity.
 Let $(C,0)$ be a germ of reduced curve on $(S,0)$, such that,  denoting  $C'$ the strict transform of $C$ on $X$, $C'$ is  a multi-germ of smooth curves on $X$.
 
 Then, using the $\Q$-cycle $Z_C$ associated  $C$ in the way defined above, 
 and denoting $D_C:=Z_C+\lfloor-Z_C \rfloor$,
 one has the following formula\footnote{Beware that in loc. cit. the $+$ before the second term in the right hand-side of the corresponding  formula (5) there is not properly printed, yet it {\em is} a plus. One should also read  formula (3) there as  $\underline{D}:= e(D_s)-\lceil D_s\rceil=e(D_s)+\lfloor - e(D_s) \rfloor$ which tallies   my definition for $D_C$.}  ~:
 $$\delta(C,0)= -\frac{1}{2} (Z_C\cdot (Z_C+Z_K))+\frac{1}{2} (D_C\cdot (D_C+Z_K)).$$
 \end{thm}
 
 Thanks to Spivakovsky's theorem~\ref{thm-mark}, we may apply the foregoing to $(C,0)\subset (S,0)$ a general polar curve of a minimal  singularity, $X$ the minimal resolution of $(S,0)$, and  $Z_C=-Z_\Omega$. As a corollary to these two theorems, we state~:
\begin{cor}
Let $(S,0)$ be a minimal singularity of normal  surface (hence rational by thm.~\ref{thm-mini-equiv-rat-reduced}).
The $\delta$ invariant of the generic polar curve is a topological invariant of $(S,0)$ i.e. depends only on the data of the weighted resolution graph.
\end{cor}

Applying  the formula in~\ref{thm-Giraud}  to get $\delta$ for the polar curve in concrete cases leads to huge computations, except in very simple~:

\begin{Example}
\label{ex-veronese-giraud}
Let  $(S,0)$ be the singularity at the vertex of the cone over a rational normal curve of degree $n$.
It is the minimal singularity whose (dual) resolution graph has only one vertex of weight $n$.
Assume that $n\geq 3$. Check that denoting $E$ the irreducible exceptional divisor, one has
$Z_\Omega=(2n-2)/n E$, $Z_K=-(n-2)/n E$, $\EZ=2E$ and hence
$\delta(C,0)=3n-6$.
\end{Example}

In section~\ref{sec-contrib-tg-cone}, we will get the result of the foregoing  example (and more)  from a geometric argument, with no use of the theorems above.
The problem of computing $\delta$ for the general polar curve of any minimal singularity is solved in~\ref{thm-scott-polar}.
















\section{A characterisation of the  generic polar curve in  a resolution}

\label{sec-char-resol}
As a consequence of the results of sections~\ref{sec-la-Note} and~\ref{sec-delta-inv},  we get the following characterisation for generic polar curves on the minimal resolution of the surface~:

\begin{thm}
\label{thm-char-resol}
Let $(S,0)$ be a minimal normal surface singularity, and $X$ the minimal resolution of  $(S,0)$.
Let $C(D)$ be {\em any  polar curve}  of $(S,0)$ with the property that its strict transform $C'(D)$ on $X$ is exactly as depicted in thm.~\ref{thm-mark}.

Then $C(D)$  is {\em a generic polar  curve} $\PSO$ as defined in def.~\ref{defi-PSO}, i.e. has the generic invariants defined in~\ref{defi-equising-data}  of the introduction.
\end{thm}

\begin{proof}
The description  of prop.~\ref{prop-descri-pol} rests only on the the shape of the polar curve in the resolution $X$, and gives in particular the datum (ii) in~\ref{defi-equising-data} (cf. def~\ref{defi-PSO} and prop.~\ref{prop-descri-discri}).
Giraud's  theorem~\ref{thm-Giraud}  gives the value of the delta invariant also from the  data of the resolution. Considering the linear system of polar curves, our special polar curve is now equisingular in the sense of def.~\ref{defi-equising-data} to the generic polar curve.
\end{proof}

\begin{rem}
\label{rem-elem-generaux-ideaux}
We explained in~\cite{gen-elem} how such characterisations of ``general" curves on a resolution may  be useful~:  here it will be used in rem.~\ref{rem-ou-on-utilise-Spivak}. 
\end{rem}





We also need  the following inductive property for  which we will use\footnote{Ideally, we would have liked not to do so, see precisely (a) of the proof of this proposition} the explicit form  of the cycle $Z_\Omega$ in (\ref{eq-ZOmega})  before Spivakovsky's thm.~\ref{thm-mark}~:

\begin{prop}
\label{prop-comport-polaire-induction}
Let $(S,0)$ be a minimal singularity of normal surface, with dual resolution graph $\Gamma$.
Let  $S_1$ be the  blow-up of $(S,0)$ at $0$ and $O_i$ a singular point of $S_1$.
Let $\Gamma_i\subset \Gamma$ be the Tyurina component corresponding to $O_i$ as in def.~\ref{defi-Tyurina-comp}. 
Let $Z_{\Omega_i}$ be the cycle associated to $\Gamma_i$ as $Z_\Omega$ is associated to $\Gamma$ in thm.~\ref{thm-mark}.

Then for any vertex $x\in \Gamma_i$  the corresponding component $L_x$ on $X$ satisfies
the following intersection property~:
\begin{equation}
\label{eq-Zomega-i}
(Z_\Omega\cdot L_x)=(Z_{\Omega_i}\cdot L_x).
\end{equation}
This means that the corresponding component $L_x$ is intersected by exactly the same number of branches of the generic polar curve for $(S,0)$ or for $(S_1,O_i)$, and the central arcs in $\Gamma_i$ are obviously also central arcs in $\Gamma$.
\end{prop}
\noindent{\it Proof.}
Although the assertion in~(\ref{eq-Zomega-i}) follows easily  from the explicit form of the cycles $Z_\Omega$ and $Z_{\Omega_i}$ (cf.~(\ref{eq-ZOmega}) p.~\pageref{eq-ZOmega}), we distinguish between~:
\begin{itemize}
 \item[(a)] the   components $L_x$ with $w(x)-v_{\Gamma_i}(x)\geq 2$.
 Since $x\in \Gamma_i$, $w(x)=v_\Gamma(x)$, hence the property in $\Gamma_i$ implies that $x$ is a central vertex in $\Gamma$. Hence $L_x$ bears components of the strict transform of the general polar curve of $(S,0)$, and here we know no other reason than computing to prove~(\ref{eq-Zomega-i}).
 
 \item[(b)] the  central components $L_x$  in $\Gamma_i$ (central vertex or  boundary of a central arc). Then, it is also central in $\Gamma$, and we believe (\ref{eq-Zomega-i})  should be understood without any reference to the cited formula, using the following remark in~\cite{Sp} p.~459 (first lines)~:
``in the neighbourhood of $L_x$, $\tilde{\Omega}$ is generated by sections whose zero set is contained in the exceptional divisor near $L_x$".\qed
 \end{itemize}

\section{The contribution of the tangent cone in the polar curve}
\label{sec-contrib-tg-cone}

In section~\ref{sec-la-Note}, we said that $\PSO$ was formed by $A_n$-curves.
Here we  explain how the bunches of $A_1$-curves arise, and  will be more precise about their geometry. 

\subsection{Discriminant and polar curve for cones over Veronese curves}
\begin{rem}
\label{rem-defi-petit-delta}
Let $(S,0)$ be the singularity of the cone over the rational normal curve of degree $m$ in $\p^m_{\C}$, whose  dual graph has just one vertex, with weight~$m$.

Denoting $P_m$ the polar curve for a generic projection of $(S,0)$ onto $(\C^2,0)$, it  is just the cone over the critical set of the projection of the rational normal curve onto $\p^1_{\C}$, which is a set of $2m-2$ distinct points by Hurwitz formula. 

Hence we know that here $P_m$ is given by $(2m-2)$ lines in $(\C^{m+1},0)$ with~:

\noindent i)  $\delta$-invariant  $3m-6$ as computed in example~\ref{ex-veronese-giraud}, from Giraud's formula.

\noindent ii) obviously a set of $2m-2$ distinct lines in $(\C^2,0)$ as generic plane projection, denoted $\delta_m$.

\end{rem}

We can say more on the geometry of $P_m$ in this case, and re-find the value of~$\delta$~:
 

\begin{lem}
\label{lem-flenner-greuel}
The general polar curve $P_m$ of the singularity  of a cone over a Veronese curve of degree $m\geq 3$, is a set of $(2m-2)$ lines in $(\C^{m+1},0)$, which 
has the generic (minimum) value of the $\delta$-invariant  for {\em any set} of $2m-2$ lines in $(\C^{m+1},0)$, and this value is  $3n-6$.
\end{lem}                                

\begin{proof}

 (a) We will  denote $V=v_m(\p^1)$ the rational normal curve of  degree $m$ in  $\p^m_C$  and $G_p(m-2,m)$ the Grassman manifold of subspaces of codimension two  in this $\p^m_C$,  and consider the map~:
$$G_p(m-2,m)\rightarrow \Hilb_V^{2m-2}$$
onto the Hilbert scheme parametrizing the  set of $2m-2$-points in $V$, which  assign to each $\Lambda$ the critical subscheme of the projection along $\Lambda$.

Using a result of H. Flenner and  M. Manaresi (in~\cite{Fl-Ma}  3.3-3.5) this map is generically finite, and since both spaces have dimension  $2m-2$ and the target space  is irreducible, the image of this map is dense.

(b) Now from a result  of  G.M. Greuel in  \cite{Gr} (3.3), a  set of $r$-lines through the origin in $\C^{m+1}$, corresponding to a set $p_1,\dots,p_r$ of points in $\p^m_\C$, has the generic $\delta$ invariant, if for all $d$ in some  bounded set of integers, their  images $v_d(p_1),\dots,v_d(p_r)$ by the corresponding Veronese embedding $v_d : \p^m_C \rightarrow \p^{N_d}_\C$ span a projective space of maximal dimension.

If we take $V\subset \p^m_\C$ to be a Veronese curve, one may always find such generic sets of points on $V$ since by composing the Veronese embeddings in Greuel's condition with the Veronese embedding defining $V$ this amounts to a genericity conditions for  points  in $\p^1_\C$.

Hence there is an open subset $U\subset \Hilb_V^{2m-2}$ with the properties that  the cone over this set of points  has the minimum delta invariant.
Applying (a) gives that these points actually occur as critical locus.

(c) A formula for the delta invariant for such a generic configuration of $r$ lines in $\C^n$ is given by Greuel in loc. cit. We leave it to the reader to check that it gives $3n-6$ in our situation.
 \end{proof}
   
\subsection{Geometry of the tangent  cone of a minimal singularity}

\begin{rem}
\label{rem-CSO-minimal}
Let $(S,0)$ be a minimal normal surface singularity with embedding dimension $N$, and $C_{S,0}$ be  its tangent cone in $(\C^N,0)$. 

Then if $\p: \C^N \setminus \{ 0\} \rightarrow \p^{N-1}_{\C}$ denotes the standard projection, the projective curve $\p(C_{S,0})$ is a non-degenerated\footnote{i.e. non contained in a proper linear subspace of $\p^{N-1}_{\C}$} curve of  minimal degree in $\p^{N-1}_{\C}$. Indeed,  condition (i) in def.~\ref{defi-kollar} immediately passes  to  $\p(C_{S,0})$.

It then follows by a standard argument (cf. e.g.~\cite{Art-deform}  p.67--68)  that each of its irreducible components is a  rational normal curve of a  linear subspace of $\p^n_{\C}$.
\end{rem}

Let $\Gamma$ be the dual graph of the minimal resolution of $(S,0)$. From Tyurina's thm.~\ref{thm-Tyurina} and the remarks following it, an irreducible component $L_{x_i}$   of $\p(C_{S,0})$ corresponds to a  vertex $x_i$ with $w(x_i)>v(x_i)$ in $\Gamma$ and it is easy to compute that the degree $m(x_i)$ of the rational normal curve $L_{x_i}$ is precisely $w(x_i)-v(x_i)$.

\begin{Conclusion}
Hence the tangent cone $C_{S,0}$ is embedded in $(\C^N,0)$  as a union of cones  over rational normal curves of degree $m_i$ intersecting along singular lines. 
\end{Conclusion}

\subsection{Scheme-theoretic polar curves and discriminants}
\label{subsec-fitting}


 To study deformations of polar curves and discriminant, we need  a scheme-theoretic definition for these objects, as introduced by B. Teissier in~\cite{Hunting} through the use of Fitting ideals.
 
 We call $P^F(S,0)$ the polar curve of a generic projection $p$ of $(S,0)$ onto $(\C^2,0)$ as defined by the Fitting ideal $F_0(\Omega_p)$ in $\OS$ and $\Delta^F_{S,0}$ its image as defined by $F_0(p_*(\Oc_{P^F(S,0)})$ in $\Oc_{\C^2,0}$.
 
  For generic projection $p$ of a normal surface $(S,0)$, these two schemes  are generically reduced so that their divisorial parts $\divis P^F(S,0)$ and $\divis \Delta^F_{S,0}$ coincide with the reduced polar curves and discriminants  defined in section~\ref{sec-pol-inv}. \footnote{But as explicitly proved in~\cite{RES} 3.5.2,  the Fitting polar curves and discriminants for minimal singularities do have embedded components as soon as the multiplicity is bigger than~$3$.}
  
  Apply these definitions to the  non-isolated singularity $(C_{S,0},0)$, we obtain~:

\begin{lem}
\label{lem-explicit-delta-CS0}
Let $(S,0)$ be a minimal normal surface singularity, with tangent cone $C_{S,0}$ ,  $\Gamma$  the dual graph of the minimal resolution of $S$. Recall that we then denote $\Gamma_{TC}$ the set of vertices $x_i$ in $\Gamma$ with $w(x_i)>v(x_i)$. 

Here, we denote by $m(x_i)$ the difference $w(x_i)-v(x_i)$, and we have just  seen that $C_{S,0}$  is made of cones over rational normal curves of degree $m(x_i)$ intersecting along singular lines.
Hence, considering the discriminant of a projection of $\C^N$ onto $\C^2$ restricted to $C_{S,0}$ we get that~:
$$\divis P^F_{C_{S,0}}=\bigcup_{x_i\in\Gamma_{TC}}  P_{m(x_i)} \cup \mbox{the singular lines in} \; C_{S,0}\, \mbox{with some multiplicity},$$
and
 $$\divis \Delta^F_{C_{S,0}}=\bigcup_{x_i\in\Gamma_{TC}} \delta_{m(x_i)} \cup \mbox{non reduced lines},$$
 \noindent where the $P_m$ and $\delta_m$  were defined in rem.~\ref{rem-defi-petit-delta} and lem.~\ref{lem-flenner-greuel}.
\end{lem}


\subsection{Deformations of the polar curves and discriminants}
\label{sub-sec-fulton}
We first  recall what we need from the construction of the deformation of  $(S,0)$ onto its tangent cone $C_{S,0}$, as described in~\cite{Fulton} chap.~5:
let $M$ be the blow-up of $(0,0)$ in $S\times \C$, and $\rho : M\rightarrow \C$ the flat map induced by composing the blow-up map with the second projection.
One then shows that~: for all $t\neq 0$ the fiber $M_t:=\rho^{-1}(t)$ is isomorphic to $S$ and $M_0$ is the sum of the two divisors on $M$, namely $\p(C_{S,0}\oplus 1)+S_1$ where $S_1$ stands for the blow-up of $S$ in $0$. To this deformation, we will apply  the following~:
\begin{prop}
\label{prop-deformation-polar-generique}
Let $\rho : X\rightarrow \D$ a flat map, with a section $\sigma$ so that  the germs  $(X_t,\sigma(t))$ are isolated singularities for $t\neq 0$, $X_0$  is a reduced possibly non-isolated singularity, and dim $X_t$ is two for all $t$.

Then, reducing the disk $\D$, one may find a projection $p : X \rightarrow \C^2\times \D$ compatible with $\rho$ so that for all $t\in \D$ the polar curve of $p_t : X_t\rightarrow \C^2\times \{ t \}$ is generic, and its image is also the  generic discriminant  $\Delta_{X_t,0}$ as defined in section~\ref{sec-pol-inv}.
\end{prop}

The trick in the proposition  above  is well-known to specialists and  may be deduced from more general results, but  we don't know an  explicit reference in the literature~: hence  a  proof will be given  in~\cite{A-B}.

Applying  the proposition to  the foregoing deformation $\rho : M  \rightarrow \D$ gives that $P^F_{S,0}$ deforms  onto $P^F_{C_{S,0}}$ and the same statement for the Fitting discriminants.
The description of the generically reduced branches of $P^F_{C_{S,0}}$ in lem.~\ref{lem-explicit-delta-CS0} now implies~:
\begin{cor}
\label{cor-descri-contrib-tg-cone}
Let $(S,0)$ be a minimal singularity, with notation as in lem.~\ref{lem-explicit-delta-CS0}, let us  denote $L_{x_i}$ the component of $\p(\CSO)$ corresponding to $x_i\in \Gamma_{NT}$. Then~:

\noindent {\rm (i)}  The generic polar curve $\PSO$ of $(S,0)$ contains a union~:
$$P_{TC}=\bigcup_{x_i\in \Gamma_{TC}} P_{m(x_i)}$$
of generic configuration of lines $P_{m(i)}$ as described in lem.~\ref{lem-flenner-greuel}.
The bunch  $P_{m(x_i)}$ in $P_{TC}\subset \PSO$ is by definition  the set of branches 
of $P_{S,0}$ which are deformed onto the (scheme-theoretically) smooth branches $P_{m(x_i)}\subset P^F_{\CSO}$ of lem.~\ref{lem-explicit-delta-CS0}.


\noindent {\rm (ii)} The same statement is true for the generic  discriminant $\DS$ of $(S,0)$~: 

 denoting
$\Delta_{TC}=\cup_{x_i\in\Gamma_{TC}} \delta_{m(x_i)},$ with $\delta_{m(x_i)}$ standing for $2m(x_i)-2$ lines in $(\C^2,0)$,  we may just as well say 
 that these smooth branches with pair-wise distinct tangents just form a $\Delta_{TC}$  part in  $\Delta_{S,0}$.

\noindent {\rm (iii)} Denote $S_1$ the blow-up of $(S,0)$.  The strict transforms on $S_1$ of the smooth curves  in  $P_{m(x_i)}\subset \PSO$   intersect the exceptional divisor only in $L_{x_i}$ and this    
intersection is transverse.


\end{cor}
\begin{proof}
(i) A curve deforming onto a smooth curve is certainly smooth, hence locally a  line. In  lem.~\ref{lem-flenner-greuel}, we said the $P_m$-curves are characterised by the minimality of $\delta$. By semi-continuity of this $\delta$ applied to the family deforming onto $P_{m(x_i)}\subset P_{\CSO}$ we get the full conclusion for the curves in $\PSO$. (ii) is direct from (i).

(iii)  Let's denote $\rho$ the deformation onto the tangent cone as recalled at the beginning of this section~\ref{sub-sec-fulton}.  The fiber $\rho^{-1}(0)$ contains the blow-up $S_1$ of $(S,0)$ intersecting $\p(\CSO\oplus 1)$ in $\p(\CSO)$.  
Since the lines  $P_{m(x_i)}$ in $\CSO$ are transverse to the Veronese curve $L_{x_i}$ in the $\p(\CSO)$ at infinity, it also follows that the strict transforms of the  curves in $P_{m(x_i)}\subset \PSO$ are transverse to the corresponding exceptional component $L_{x_i}\subset \p(\CSO)$ on the blow-up $S_1$.
\end{proof}

\section{Limit trees}
\label{sec-limit-tree}

We  proceed to identify the remaining part in $\PSO$ besides the $P_{TC}$-part just exhibited. The following  
{\em limit tree} construction introduced by  T. de Jong and D. van Straten in~\cite{DJ-VS} will turn out to be much relevant to this description.
Precisely,   using the  height function we defined in~\ref{defi-height}, one finds in loc. cit. (1.13) the~:

\begin{defi}
\label{defi-limit-equiv}
Let $\Gamma$ be the dual graph of minimal resolution of a minimal singularity of normal surface.
A {\em limit equivalence relation} $\sim$ is an equivalence relation on the vertices of $\Gamma$ satisfying the following two conditions~:

(a) Vertices $x$ with height $s_x=1$ i.e. with $w(x)>v(x)$ belong to different equivalence classes,

(b) for every vertex $x$ in $\Gamma$ with height $s_x=k+1$, $k\geq 1$, there is exactly {\em one vertex} $y$ connected to $x$ with height $s_y=k$ and $y\sim x$.

Then, the tree $T=\Gamma/\sim$ is a called a {\em limit tree} associated to $\Gamma$.

It is clear that any equivalence class contains exactly one vertex $x_i$ of height one, so that we denote this equivalences classes as vertices $\widetilde{x_i}$ in $T$.
 \end{defi}

In fact, we only make this construction in the particular case of minimal singularities with {\em reduced graphs}  in the sense of notation~\ref{notation-graphe}, so that the definition above really correspond to the definition in loc. cit.\footnote{For the non-reduced case, one has to use the extended resolution graph of loc. cit. to build the limit tree, to really get a bijection between vertices of $T$ and element of the set $\mathcal{H}$ considered in loc. cit. But, again, we won't use this.}

Starting with $\Gamma$ as in example~\ref{Example-montre-height} , one may associate non-isomorphic limit trees to the same reduced graph $\Gamma$, depending on the equivalence classes chosen, namely~:
\begin{figure}[h]
\caption{\label{Two-limit-trees} Two distinct  limit trees for the dual graph in Figure~\ref{Graphe-a-deux-limites}.}
\unitlength=1cm

\begin{picture}(5,2)

\put(0.01,1){$T_{1}:$}

\put(1,1){$\ast$}

\put(1,1.3){$\tilde{x}_{1}$}

\put(1.1,1.1){\line(1,0){1}}

\put(2.1,1.){$\ast$}

\put(2,1.3){$\tilde{x}_{4}$}

\put(2,1.1){\line(1,0){1}}

\put(3,1){$\ast$}

\put(3,1.3){$\tilde{x}_{2}$}

\put(3,1.1){\line(1,0){1}}

\put(4,1){$\ast$}

\put(4,1.3){$\tilde{x}_{3}$}

\end{picture}


\unitlength=1cm

\begin{picture}(5,2)

\put(0.1,1){$T_{2}:$}

\put(1,1){$\ast$}

\put(1.1,1.1){\line(1,0){1.1}}

\put(1,1.3){$\tilde{x}_{1}$}

\put(2,1){$\ast$}

\put(2,1.3){$\tilde{x}_{2}$}

\put(2.1,1.1){\line(1,0){1}}

\put(3,1){$\ast$}

\put(3,1.3){$\tilde{x}_{4}$}

\put(2.1,1){\line(0,-1){1}}

\put(2,0){$\ast$}

\put(2.3,0){$\tilde{x}_{3}$}

\end{picture}\\          
\end{figure}
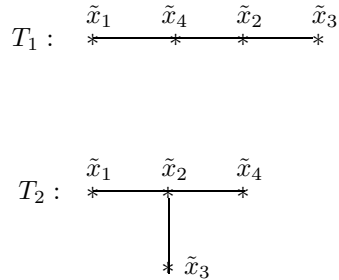

\begin{Notation}
For any pair $x,y$ of vertices on the dual graph $\Gamma$, we denote by $C(x,y)$ the (minimal) chain on $\Gamma$ joigning them (including the end points). (It is unique since $\Gamma$ is  a tree).

We define the length $l(x,y)$ to be the number of vertices on $C(x,y)$ and the overlap $\rho(x,y;z)$ as the number of vertices on $C(x,z)\cap C(y,z)$.
\end{Notation}

As in \cite{DJ-VS}, we attach to a limit tree $T$  the  following data~:
\begin {itemize}
\item for any edge $(\tilde{x},\tilde{y})$ of  $T$, the length $l(x,y)$ where $x,y$ are the corresponding vertices of heigth one in the resolution  graph $\Gamma$,
\item for any pair of adjacent edges  $(\tilde{x},\tilde{z})$ and $(\tilde{z},\tilde{y})$ in $T$, the overlap $\rho(x,y;z)$.
\end{itemize}

We use the notation $(T,l,\rho)$ for the data above. In loc. cit. lemma (1.16), it is shown that these data determine uniquely the resolution graph $\Gamma$.
 

\section{Description of the polar curve using the limit tree}
\label{sec-discri-limit-tree}





The following is our main result, we formulate it for the polar curve $\PSO$  reminding  that this implies the analogous statements for the discriminant $\DS$~:

\begin{thm}
\label{main-thm}
Let $(S,0)$ be a minimal singularity of normal surface. Let $\Gamma$ be the dual graph of the minimal resolution of $S$.

Let $\Gamma^r$ be  the reduced graph associated to $\Gamma$ in the sense of notation~\ref{notation-graphe} i.e. the same graph with the weights of the $x_i$ of height one reduced to $v(x_i)+1$, and let $(S^r,0)$ be a minimal singularity with  dual resolution graph $\Gamma^r$.

Then  the generic polar curve  $\PSO$  decomposes into~:
$$\PSO=P_{TC}\cup P_{S^r}$$
where the contact between any line in $P_{TC}$ and any branch in $P_{S^r,0}$ is one
and $P_{TC}$ was described in cor.~\ref{cor-descri-contrib-tg-cone} as the ``contribution of the tangent cone''.


Let $T$ be  limit tree for $\Gamma^r$, as defined in section~\ref{sec-limit-tree} and $(T,l,\rho)$ the set of data (length and overlap) associated  to it at the end of that  section.

These data give the following  easy description of  $P_{S^r}$ (as a union of $A_n$-curves)~:
\begin{trivlist}
\item{$\bullet$} each edge $(\tilde{x_i},\tilde{x_j})$ in the limit tree $T$ defines exactly one   $A_{l_{i,j}}$-curve in $P_{S^r}$, where $l_{i,j}$ stands for $l(\tilde{x_i},\tilde{x_j})$. 
\item{$\bullet$} For each pair of adjacent edges $(\tilde{x_i},\tilde{x_j})$ and $(\tilde{x_j},\tilde{x_k})$, the contact (def.~\ref{defi-contact}) between the corresponding $A_{l_{i,j}}$ and $A_{l_{j,k}}$-curves in  $P_{S^r}$ is exactly the overlap $\rho(i,k;j)$.
\item{$\bullet$} For non adjacent edges $(\tilde{x_i},\tilde{x_j})$ and $(\tilde{x_k},\tilde{x_l})$, the contact between the corresponding $A_{l_{i,j}}$ and $A_{l_{k,l}}$-curves in $P_{S^r}$  is the {\em minimum} of the contacts between adjacent edges on the chain joining them.
\end{trivlist}
\end{thm}

Let us first illustrate this on the following~:
\begin{Example}
(i) For  a minimal singularity $(S,0)$ with dual graph as in Figure~\ref{Graphe-a-deux-limites} p.~\pageref{Graphe-a-deux-limites}, using any of the limit trees in Figure~\ref{Two-limit-trees}, we get~:
$\PSO=A_5\cup A_5'  \cup A_3$, with contact three between the two $A_5$ and contact one between the $A_5$'s and the $A_3$.

\noindent (ii)  For the example~\ref{Ex-la-Note}, the description of the discriminant was already given there.
  It is now more directly seen from the limit tree in Figure~\ref{fig-limit-ex-note} given below together with the data $(l,\rho)$, where the lengths $l$ are put above the edges and the $\rho$ as smaller numbers in-between a pair of edges  (following the same convention as in~\cite{DJ-VS} (1.19)).
\end{Example}
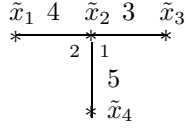
\begin{figure}[h]
\caption{Limit tree for the reduced  graph associated to the graph in Example~\ref{Ex-la-Note} \label{fig-limit-ex-note}}

\unitlength=1cm

\begin{picture}(5,2)


\put(1,1){$\ast$}

\put(1.1,1.1){\line(1,0){1.1}}

\put(1.5,1.3){$4$}

\put(1,1.3){$\tilde{x}_{1}$}

\put(2,1){$\ast$}

\put(1.8,0.8){$\scriptstyle 2$}
\put(2,1.3){$\tilde{x}_{2}$}
\put(2.2,0.8){$\scriptstyle 1$}

\put(2.1,1.1){\line(1,0){1}}

\put(2.5,1.3){$3$}

\put(3,1){$\ast$}

\put(3,1.3){$\tilde{x}_{3}$}

\put(2.1,1){\line(0,-1){1}}

\put(2.3,0.4){$5$}

\put(2,0){$\ast$}

\put(2.3,0){$\tilde{x}_{4}$}

\end{picture}\\

\end{figure}

\noindent {\em The rest of this section is devoted to the proof of thm.~\ref{main-thm} above}.

First we recall the following well-known~:

\begin{lem}
\label{lem-mult-discri}
Let $(S,0)$ be a minimal singularity of normal surface and $m$ be the multiplicity of $(S,0)$.
Then the multiplicity of the generic polar curve (resp. discriminant) is $2m-2$
\end{lem}

\begin{proof}
This is easily deduced from the two following facts (we refer e.g.  to~\cite{RES} (3.9) and  \S~5)~:

\noindent (a) for any normal surface $(S,0)$ and  any projection $p : S\rightarrow \C^2$ whose degree equals the the multiplicity $m$ of the surface, the  multiplicity of the discriminant $\Delta_p$ is $m+\mu-1$ where $\mu$ is the Milnor number of a generic hyperplane section of $(S,0)$.

\noindent (b) When $(S,0)$ is minimal, then  $\mu=m-1$.
\end{proof}

\noindent The proof of thm.~\ref{main-thm}  is by induction on the maximal height of the vertices in $\Gamma$~:

\medskip
{\bf A) Initial step~--}  The maximal height in $\Gamma$ is one.
  We prove the result by a direct argument (independent from Spivakovsky's theorem).
Now,  all the vertices $x_i$ in $\Gamma$  are in $\Gamma_{TC}$,
and the minimal resolution $X$ of $(S,0)$ is the first blow-up.

\noindent (a )We know from the deformation onto the tangent cone that each exceptional component $E_{x_i}$  bears the strict transform of $(2n_i-2)$ smooth branches of the polar curve, cutting $E_{x_i}$ transversely at general points,  with $n_i=w_i-v_i$ (cf.~\ref{cor-descri-contrib-tg-cone} (iii)).

\noindent (b) A  general theorem of J.  Snoussi (in \cite{Sn}), valid for any normal surface singularity, describes the base points of the linear system of polar curves on the first normalized blow-up of $(S,0)$.   In our situation the blow-up is already normal and even  smooth, and hence Snoussi's theorem implies that here the   bases points are exactly the singular points of the exceptional divisor, i.e. the intersection points of two components $E_{x_i}$ and $E_{x_j}$. 

Let $N$  be the number of vertices in $\Gamma$, then $\Gamma$ has $N-1$ edges (it is a tree), which  represent the intersections points of exceptional components.

By Bertini's theorem,  the part of the generic polar curve $\PSO$ whose strict transform goes through a base point  is singular, i.e. has multiplicity at least two.

Hence, adding the contributions  of the smooth branches in (a)  and the singular curves in (b), 
the multiplicity $m(\PSO,0)$ of the polar curve satisfies the  inequality~:
\begin{equation}
\label{ineq-eq}
m(\PSO,0)\geq \sum_{i=1}^N (2n_i-2)+ 2(N-1).
\end{equation}
Comparing this to the equality $m(\PSO,0)=2m-2$ of  lemma~\ref{lem-mult-discri} above, where the multiplicity $m$ of $(S,0)$ equals the $\sum_{i=1}^N n_i$,  proves that (\ref{ineq-eq}) is in fact in equality.

Hence, each point of intersection of two exceptional components bears  a curve of multiplicity exactly two. We now claim that the curve in question is the strict transform of  a $A_2$-curve singularity on $(S,0)$.
Let  $C$ be such a curve.

Then, the multiplicity  $m(C,0)=2$ is the intersection number of $C$ with a generic hyperplane section of $(S,0)$. This intersection number may be computed on $X$ as the intersection number of the strict transform $C'$ with the reduced exceptional divisor (which is the cycle defined by the  maximal ideal of $(S,0)$). Since we know $C'$ intersects two exceptional components, the intersection of $C'$ with each one should be transverse.

Hence $C$ is a branch a multiplicity two resolved in one blow-up, i.e. $A_2$-curve.
This completes the proof of the initial step A.


\medskip

\noindent {\bf B) The induction step~--}  We use first the following general lemma in~\cite{RES} 6.1~:
 
\begin{lem}
\label{lem-RES}
Let $(S,0)$ be any normal surface singularity and $p : (S,0)\rightarrow (\C^2,0)$ any projection with degree equal to  the multiplicity $\nu=m(S,0)$.

Then, denoting $b_0 : \widetilde{\C^2}\rightarrow (\C^2,0)$ the blow-up of the origin, and $\Sigma_1$ the analytic fiber product of $b_0$ and $p$ above $(\C^2,0)$, one proves that the normalisation of $\Sigma_1$ coincides with the normalized blow-up $\overline{S_1}$ of $(S,0)$, which yields   the following commutative diagram~:

$$\xymatrix{\overline{S_1}\ar@/_/[ddr]_{\phi_1} \ar[dr]^n \ar@/^/[drr]^{\overline{b_1}} && \\&\Sigma_1\ar [d]^{p_1} \ar[r] & \ar[d]^p (S,0)  \\ &\widetilde{\C^2} \ar[r]^{b_0} & (\C^2,0)}$$

\noindent where $\varphi_1 : \overline{S_1}\rightarrow  \widetilde{\C^2}$ is  the composition of the pulled-back projection $p_1$ with the normalization $n$.\par

The following formula can then be obtained for the discriminant  $\Delta_{\varphi_1}$~:
\begin{equation}
\label{eq-RES}
\Delta_{\varphi_1}=(\Delta_p)'+(\nu-r) E,
\end{equation}
where $(\Delta_p)'$ is the strict transform of the discriminant of $p$,  $E$ denotes the reduced exceptional divisor, $\nu$ is the multiplicity of the germ $(S,0)$ and $r$ the number of branches of a general hyperplane section of $(S,0)$.\par
\end{lem}

We refer to loc. cit. for the proof, we just precise that the discriminants in  the equality of the lemma are the divisorial parts of Fitting discriminants as defined in section~\ref{subsec-fitting}, which are hence allowed to have non-reduced components.

Here, $(S,0)$ being a minimal singularity, the blow-up $S_1$ is already normal (cf. e.g. \cite{RES} Thm. 5.9), so that  $\overline{S_1}$ is just $S_1$. The generic projection that we consider fulfils certainly the property $\dgre p=m(S,0)$ as a necessary condition. 
Since the general hyperplane section of a minimal singularity of multiplicity $\nu$ is just $\nu$ lines (cf. e.g. loc. cit. lem.~5.4), formula~(\ref{eq-RES})  in the above lemma simply reads~: 
\begin{eqnarray*}
\label{eq-ST}
\Delta_{\varphi_1}=(\Delta_p)',
\end{eqnarray*}
and similarly, denoting $C(D)$ the  polar curve of the projection $p$, $C'(D)$ its strict transform on $S_1$  and $C_{\varphi_1}$ the polar curve for the projection $\varphi_1$, we get~:
\begin{equation}
\label{eq-ts}
C'(D)=C_{\varphi_1}.
\end{equation}

 From thm.~\ref{thm-Tyurina} (see also def.~\ref{defi-Tyurina-comp}), we know that the singularities  $O_i$ of $S_1$ are minimal singularities whose resolution graphs are the Tyurina components $\Gamma_i$.
 
 Localising  the  result  of (\ref{eq-ts}) in  $O_i$ yields the  following~:
 
 \begin{Conclusion}
 \label{conc-ts}
 Let $C(D)$ be a generic polar curve for $(S,0)$ and $C'(D)$ its strict transform on the  blow-up $S_1$
 of $S$ at $0$. Let $O_i$ a singular point  of $S_1$.
  We proved that  the part  of $C(D)'$ going through $O_i$ is the  polar curve for the projection $\varphi_1$ obtained of the germ $(S_1,O_i)$ onto a plane, as in the lemma above.
 \end{Conclusion}
 
 \begin{rem}
\label{rem-ou-on-utilise-Spivak}
 To apply induction, we need to know that the projection $\varphi_1 :  (S_1,O_i) \rightarrow (\C^2,0)$  in question is generic, i.e. has the generic polar curve.
 
 Counting multiplicity as in A) gives that this  projection has degree equals the multiplicity of $(S_1,O_i)$, but this is not  enough to prove that the polar curve is generic.
This will be proved thanks to the results of section~\ref{sec-char-resol}.
 \end{rem}

Indeed, once we know from conclusion~\ref{conc-ts} that  $C'(D)$ is a polar curve for  $(S_1,O_i)$ we may use prop.~\ref{prop-comport-polaire-induction} to see   that the strict transform of $C'(D)$ on $X$ which is also part of the strict transform of $C(D)$, actually fulfils the conditions  of  the characterisation in thm.~\ref{thm-char-resol}. Then~:

\begin{Conclusion}
\label{conc-ts-gene}
With the same notation as in conclusion~\ref{conc-ts}, the part of $C'(D)$ going through $O_i$ is the {\em  generic polar curve} $P_{S_1,O_i}$ for the germ $(S_1,O_i)$.

\end{Conclusion}

Now, the induction hypothesis applied to each $(S_1,0_i)$ yields that 
 $P_{S_1,0_i}$ is a union of $A_{n}$-curves described by  a limit tree $T_i$ for $\Gamma_i$ as stated in  theorem~\ref{main-thm}.
 
 \medskip
 
{\noindent \bf C) Reconstructing $P_{S,0}$ from its strict transform}

\smallskip

Let  $(S,0)$ be a minimal surface singularity and let $S_1,$ be its  blow-up, and $E$ the exceptional divisor  with components $E_1,\dots,E_r$. We will denote by $O_1,\dots,O_s$ the  singular points of $S_1$ and $Q_1,\dots,Q_t$ the points of intersection of components of $E$ which are not  singular points of $S_1$.

We already know that the generic polar curve $\PSO$  of $(S,0)$ is precisely made of~:

\begin{enumerate}
\item $A_1$-curves in number $\sum_{i=1}^r (m_i-1)$, whose strict transforms intersect each $E_i$ as $(2m_i-2)$ lines going through generic points of $E_i$, for $i=1,\dots,r$,
\item $A_2$-curves singularities  in number $t$, each one having its strict transform on $S_1$  intersecting a different point $Q_i$ defined above,
\item curves whose strict transforms go through the  singular points $O_i$ of $S_1$.
\end{enumerate}

The first two points are proved by the same reasoning as in step A).
Step B) applied to curves in (3) for each $O_i$  gives the description of the strict tranforms of these curves as $A_n$-curves described by the limit tree $T_i$ associated to $(S_1,O_i)$.  

The corresponding description, for all the  curves in (3) whose strict transform go through the same $O_i$,  on $(S,0)$ itself  is then obtained by adding $2$ to all the $n$'s and one  to the $\rho$ by elementary properties of these $A_n$-curves and our def.~\ref{defi-contact} of the contact.




  But now from \cite{DJ-VS} (1.18) we know that the data associated to limit tree  $T_i$  of $\Gamma_i$ are related to $T$ exactly the same way (length:= length$-2$, overlap := overlap $-1$).
  
This completes the proof by induction for the first two points of theorem~\ref{main-thm}, the last point follows by definition of the contact. 






 
 \section{Scott deformations and polar invariants}
 \label{sec-scott}
 
 The following was first proved by de Jong and  van  Straten  in~\cite{DJ-VS} Thm.~2.13~:
 
 \begin{thm}
 \label{thm-scott}
 Let $(S,0)$ be a minimal singularity of normal surface with multiplicity $m$.
 Let $S_1$ be the blow-up of $0$ in $S$, with singular points  $O_1,\dots,O_r$.
 Then there exists a one-parameter deformation $\rho : X\rightarrow \D $ of $(S,0)$ on the Artin component such that
$X_s$  for $s\neq 0$ has $r+1$ singular points isomorphic respectively to the $(S_1,O_i)$ for $i=1,\dots,r$ and to the cone over the rational normal curve of degree~$m$.
\end{thm}
 

This has to be compared to the a standard result for plane curves, attributed to C. A. Scott in~\cite{Sand-DJ-VS}, where a proof is also given (see p.~460)~:

 \begin{lem}
 \label{lem-scott}
 Let $(C,0)\in (\C^2,0)$ be a plane curve singularity of multiplicity $m$. 
 Let $O_i$ for $i=1,\dots,r$ be the singularities of the first blow-up $C_1$  of $(C,0)$
 There exists a one-parameter  $\delta$-constant deformation  $\Gamma$ of $(C,0)$ such that $\Gamma_s$ for $s\neq 0$ is a plane curve which has $r+1$ singular points isomorphic respectively to the $(C_1,O_i)$
 for $i=1,\dots,r$ and to an ordinary $m$-tuple point. 
 \end{lem}

 Beyond the formal analogy between thm.~\ref{thm-scott} and lem.~\ref{lem-scott}, de Jong and van Straten prove the result in thm.~\ref{thm-scott}  for the more general class of {\em sandwiched singularities} as a consequence of their theory of {\em decorated curves}~:  all the deformations of these surface singularities can be obtained from  deformation of {\em decorated curves} associated to the singularity. In particular,  the Scott deformation of a decorated curve (conveniently adjusted) gives rise to the deformation  in~thm.~\ref{thm-scott}.

 As an  application of our description for generic discriminants in thm.~\ref{main-thm}, however, we get first a new relation between these two deformations~:
 
 \begin{cor}
 \label{cor-scott-discri}
 Let  the  notation be the same as  in thm.~\ref{thm-scott}. We  will also call the deformation $X$
 the  {\em Scott deformation} of the surface $(S,0)$.
 
 Considering  a projection $p$ of  $X$ in $\D\times ÷C^2$ as in prop.~\ref{prop-deformation-polar-generique} i.e. compatible with~$\rho$ and such that the discriminant $\Delta(p_t : X_t \rightarrow \C^2)$ is the generic discriminant $\Delta_t$ for all the singularities in $X_t$, for all $t$, one gets a  deformation
 $\rho' : \Delta  \rightarrow \D$ of the generic discriminant  $\Delta_{S,0} $ of $(S,0)$ which is exactly the Scott deformation of this curve as defined in lem.~\ref{lem-scott}.
 \end{cor}
 \begin{proof}
In our proof in section~\ref{sec-discri-limit-tree},  it is proved that the discriminant  of $(S_1,O_i)$ is the part of the strict transform  of  the discriminant of $(S,0)$  going through  the image of $O_i$ in the plane (it is of course also obvious from the result there).

 The discriminant of the cone over the $m$-th Veronese curve is  a  $2m-2$-tuple ordinary point in the plane (cf. rem.~\ref{rem-defi-petit-delta}).  This is indeed the last singularity occuring in the Scott deformation of $\DS$, since, by lem. ~\ref{lem-mult-discri},
the multiplicity of the $\DS$ is  $2m-2$.
 \end{proof}

Considering polar curves in this Scott deformation, we get the more interesting~:

\begin{thm}
\label{thm-scott-polar}
Let the notation be as in cor.~\ref{cor-scott-discri}. Then, the  polar curve for $p_t : X_t \rightarrow (\C^2,0)$ is also the generic polar curve $P_{X_t}$ (which is a multi-germ of space curves for $t\neq 0$).
Further,  $P_{X_t}$ is a $\delta$-constant deformation of the generic polar curve $P_{S,0}$.
Hence iterating Scott deformations, one may compute the $\delta$-invariant of $P_{S,0}$ as sum of $\delta$-invariants for sets of generic lines $P_m$ as in lem.~\ref{lem-flenner-greuel}.
 \end{thm}
\begin{proof}
In theorem~\ref{thm-scott}, the deformation $X_t$ is said to belong  to the  Artin-component of $(S,0)$.
This means that it has a simultaneous resolution, in which (cf.  lem~\ref{lem-res-cp}) the $P_{X_t}$ are also resolved. One then has a {\em normalisation in family} for the family $P_{X_t}$, which is equivalent to ``$\delta$-constant'' (cf.~\cite{Hunting}  p.~609).
\end{proof}

\begin{figure}[h]
\caption{Graph with weights on the vertices for example~\ref{ex-join-A5} \label{fig-join-A5}}

\begin{center}
\unitlength=1cm

\begin{picture}(6,8)

\put(1.1,5.1){\circle*{.1}} 

\put(1,4.7){2}

\put(1.1,5.1){\line(1,0){1}}

\put(2.1,5.1){\circle*{.1}} 


\put(2.1,4.7){2}

\put(2.1,5.1){\line(1,0){1}}

\put(3.1,5.1){\circle*{.1}}


\put(2.8,4.7){3}

\put(3.1,5.1){\line(0,-1){1}}

\put(3.1,4.1){\circle*{.1}}

\put(2.8,4.1){2}


\put(3.1,4.1){\line(0,-1){1}}

\put(3.1,3.1){\circle*{.1}} 

\put(2.8,3){2}







\put(3.1,5.1){\line(1,0){1}}

\put(4.1,5.1){\circle*{.1}}

\put(4.1,4.7){2}

\put(4.1,5.1){\line(1,0){1}}

\put(5.1,5.1){\circle*{.1}} 

\put(5,4.7){2}

\end{picture}
\end{center}
\end{figure}
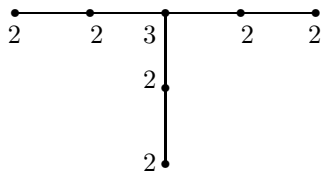

We illustrate the second statement in thm~\ref{thm-scott-polar} by giving  an~:

\begin{Example}
\label{ex-join-A5}
Taking  a singularity with graph as in figure~\ref{fig-join-A5}, and 
applying twice the Scott deformation of the surface, one gets two cones over a cubic and two cones over a conic. Hence the polar curve deforms onto two $P_3$ and two $P_2$  (in the notation of lem.~\ref{lem-flenner-greuel}), which gives $8$  for the $\delta$-invariant.\footnote{Beware that  $\delta(P_2)=1$ is not given by the formula $\delta(P_n)=3n-6$, valid for $n\geq 3
$.}
\end{Example}

Let us  end this with the following~:
\begin{rem}
\label{rem-deform}
The information on (the resolution graph of) $(S,0)$ given by the generic discriminant $\DS$ is of course partial~: e.g. one may permute the Tyurina components in the resolution graph of $(S,0)$ or the weights on the tangent cone without changing $\DS$. However, when one looks at deformations on $(S,0)$, we believe the information on the discriminant is most valuable~:

(a)  As a very basic  occurence of this~:  a  family of normal surfaces $S_t$  with constant generic discriminant $\Delta_{S_t,0}$ is  Whitney-equisingular  and in particular has constant topological type (encoded by the minimal resolution graph).
As a consequence of our result, these three equisingularity conditions are in fact equivalent for minimal singularities of surfaces (cf.~\cite{A-B}).



(b) Much more generally, one can deform the discriminant $\DS$  and ask which deformation of $(S,0)$  ``lies above" the curve-deformation~: for example, can one deduce  the existence of the Scott deformation of the  surface $(S,0)$ in the sense of thm.~\ref{thm-scott}  as deformation ``lying above" the Scott deformation of $\DS$ ?

This would give a description of some deformation theory of the surface through an invariant which,
 as opposed to the birational join construction of Spivakovsky or the decorated tree construction of De Jong and Van Straten, is uniquely defined from $(S,0)$.

\end{rem}

\end{document}